\documentclass[10pt]{article}

\usepackage{amsfonts}
\usepackage{amssymb}
\usepackage{amsmath}
\usepackage{url}
\usepackage{pdfpages}
\usepackage{indentfirst}
\usepackage{color}

\numberwithin{equation}{section}
\numberwithin{equation}{subsection}

\numberwithin{figure}{section}
\numberwithin{figure}{subsection}

\definecolor{myblue}{rgb}{0,0,1}

\newcommand{\abs}[1]{|#1|}

\newcommand{\w}{{\mathrm w}}
\newcommand{\W}{{\mathrm W}}

\newcommand{\R}{{\mathrm{R}}}
\newcommand{\RH}{{\mathrm{RH}}}

\newcommand{\myd}{{\mathrm d}}

\newcommand{\lS}{{\mathrm{Spec}}^\lambda}
\newcommand{\mS}{{\mathrm{Spec}}^\mu}
\newcommand{\lamS}{{\mathrm{Spec}}^{\ln|{\mu}|}}
\newcommand{\lamSlow}{{\mathrm{Spec}}^{\ln|{\mu}|<}}
\newcommand{\lamSup}{{\mathrm{Spec}}^{\ln|{\mu}|>}}

\newcommand{\zetastar}{{\zeta^*}}

\author{{\sc Yu.\,V.\,Matiyasevich}}

\title{{
\bf Hidden Life of \\Riemann's Zeta Function}\\
 2. Electrons and Trains}

\date{Steklov Institute of Mathematics at St.Petersburg, Russia\\[2mm]
\url{http://logic.pdmi.ras.ru/~yumat}}

\begin{document}

\maketitle

\begin{abstract}
The Riemann Hypothesis can be reformulated
as statements about the eigenvalues of
certain matrices whose entries are defined
in terms of the Taylor coefficients of the
zeta function. These eigenvalues exhibit
interesting visual patterns allowing one to
state a number of conjectures.

The Hankel matrices introduced here are
obtained, by rearranging of columns, from
Toeplitz matrices whose eigenvalues were
considered in
\url{http://arXiv.org/abs/0707.1983}. The
present paper is a continuation of this
paper and references such as 1X, 1.x,
(1.x.y), and 1.x.y refer respectively to
conjectures, sections, formulas, or figures
from there.

\end{abstract}

\section*{}
\setcounter{section}{2}
\subsection{Hankel matrix representation}

As it was indicated in Section 1.6,
 we can give an alternative representation for the
numbers $\tau_m$ (defined by the expansion (1.5.1)),
  in the form of  determinants, by
rearranging the columns of the matrices
$L_{1,m}$  introduced by  (1.6.3). More
generally, we rearrange the columns of the
matrices $L_{l,m}(f)$  defined by (1.14.2):
\begin{eqnarray}
\lefteqn{M_{l,m}(f)=}\nonumber\\
&&-(-1)^{\frac{(m+1)(m+2)}{2}}
\begin{pmatrix}
\theta_{l+m-1}(f)&\theta_{l+m-2}(f)&\dots&\theta_{l}(f)\\
\theta_{l+m-2}(f)&\theta_{l+m-3}(f)&\dots&\theta_{l-1}(f)\\
\vdots&\vdots&\ddots&\vdots\\
\theta_{l}(f)&\theta_{l-1}(f)&\dots&\theta_{l-m+1}(f)
\end{pmatrix}.
\end{eqnarray}
Clearly,  for all $l$ and $m$ we have
\begin{equation}
\det(M_{l,m}(f))=\det(L_{l,m}(f)),
\end{equation}
so we have the following reformulations of
versions 2 and 3 of subhypothesis $\RH_l^\w$
from Sections 1.14 (with $\W_l$ defined by
(1.12.10)):

\

{\bf  $\mathbf{RH_{\emph{l}}^w}$ (version
$\mathbf{2}'$)}.  \emph{For $m\rightarrow
\infty$
\begin{equation}
\det(M _{l,m}(\zetastar))\,
=\,\W_l^m(\R_l(\zetastar)+o(1)).
\label{detmulR1}\end{equation} with some
constant $\R_l(\zetastar)$}.

\

{\bf  $\mathbf{RH_{\emph{l}}^w}$ (version
$\mathbf{3}'$)}.  \emph{
\begin{equation}
\lim_{m\rightarrow
\infty}\left(\det(M_{l,m}(\zetastar))\right)^\frac{1}{m}\,
=\,\W_l. \label{detmuulRl}\end{equation} }

\subsection{Yet More Eigenvalues}

By analogy with (1.15.1), we have the
representation
\begin{equation}
\det(M_{l,m}(f))=\mu_{l,m,1}(f)\mu_{l,m,2}(f)\dots
\mu_{l,m,m}(f) \label{lmus}\end{equation}
where $\mu_{l,m,1}(f), \ \mu_{l,m,2}(f),\
\dots,
 \ \mu_{l,m,m}(f)$
are the eigenvalues of the matrix
$M_{l,m}(f)$. Respectively, we have

\

{\bf  $\mathbf{RH_{\emph{l}}^w}$ (version
5)}.
\begin{equation}
\lim_{m\rightarrow \infty}
\left(\prod_{n=1}^m{\mu_{l,m,n}(\zetastar)}\right)^\frac{1}{m}=
\W_l. \label{vers5}\end{equation}

\

The (multi)set $\{\mu_{l,m,1}(f),\
\mu_{l,m,2}(f),\dots \mu_{l,m,m}(f)\}$ will
be called $\mu$-\emph{spectrum} of the
function $f$ and will be
denoted~$\mS_{l,m}(f)$.

\subsection{Positions of the $\mu$
eigenvalues}

According to \eqref{vers5}, the (geometric)
mean of $\mu_{1,m,1}(\zetastar)$,
$\mu_{1,m,2}(\zetastar)$, \dots,
$\mu_{1,m,m}(\zetastar)$ approaches $\W_l$
when $m$ goes to infinity, which is similar
to the behavior of the eigenvalues
$\lambda_{1,m,1}(\zetastar)$,
$\lambda_{1,m,2}(\zetastar)$, \dots,
$\lambda_{1,m,m}(\zetastar)$. However, there
are many differences between the
distribution of the eigenvalues from spectra
$\mS_{l,m}(\zetastar)$ and
$\lS_{l,m}(\zetastar)$.

The first such difference is evident: the
numbers $\mu_{1,m,1}(\zetastar)$,
$\mu_{1,m,2}(\zetastar)$, \dots,
$\mu_{1,m,m}(\zetastar)$, being the
eigenvalues of a Hankel matrix with real
entries, are all real themselves.

Computations suggest that in contrast to the
case of the $\lambda$-spectra, the union
$\cup_{m=1}^{\infty}\mS_{l,m}(\zetastar)$ is
bounded neither from above nor from below.
Moreover, the point $0$ is a limit point of
this set. That is why it is  reasonable to
consider the sets
\begin{equation}
\lamS_{l,m}(f)=\{\ln\abs{\mu}:\mu \in
\mS_{l,m}(f)\}
\end{equation}
which will be called \emph{logarithmic
$\mu$-spectra}. When exhibiting several
logarithmic $\mu$-spectra, we will shift the
lines vertically,  that is,
an eigenvalue $\ln\abs{\mu}$ from $\lamS_{l,m}(f)$
will be placed at point
$(x,y)=(\ln\abs{\mu},m)$.

\begin{figure}[phtb]
\includegraphics[width=1\textwidth]{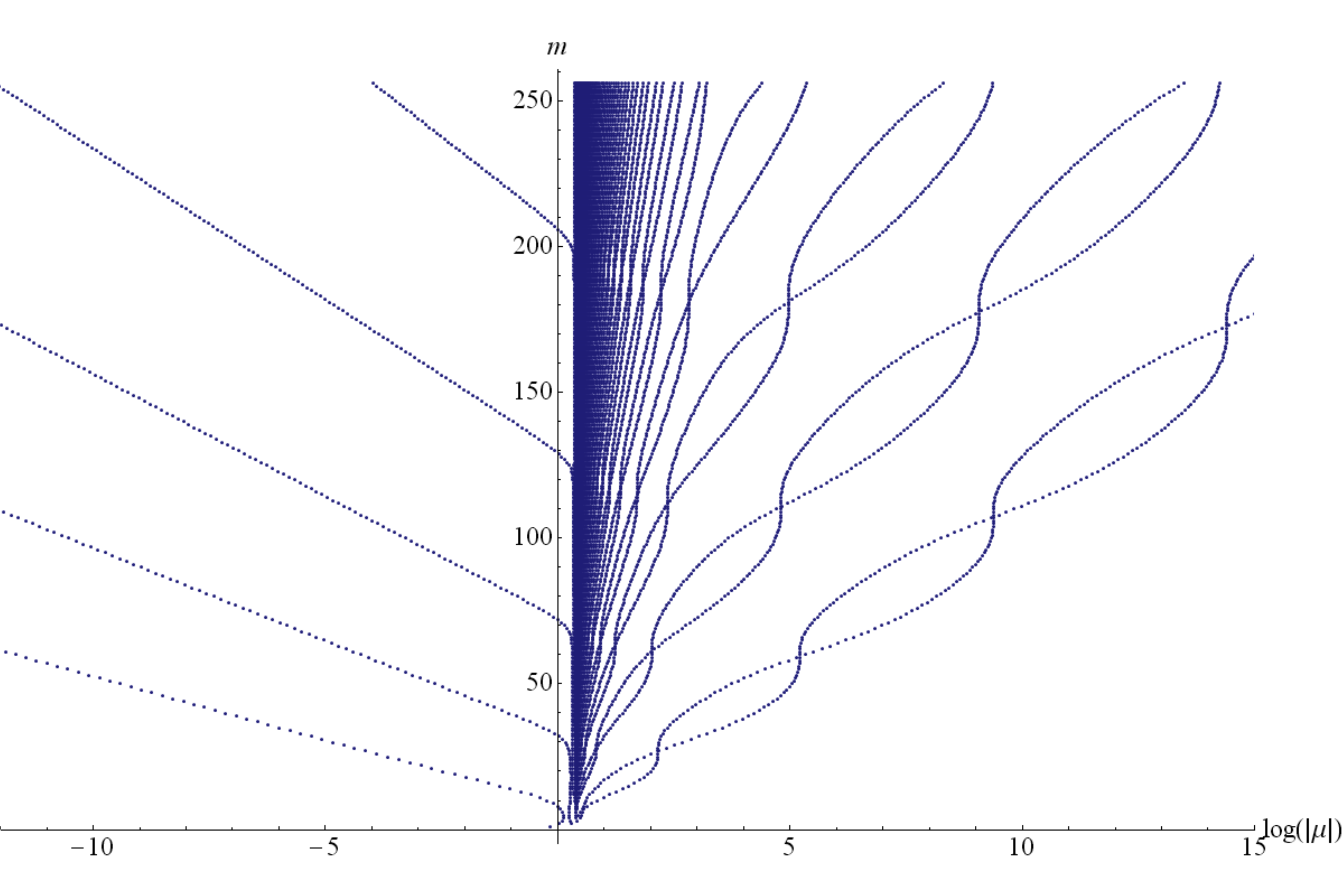}\\
\parbox{.94\textwidth}{\caption{$\lamS_{1,m}(\zetastar),\ m=1,\dots,256.$}
\label{picm1}}
\end{figure}
\begin{figure}[phtb]
\includegraphics[width=1\textwidth]{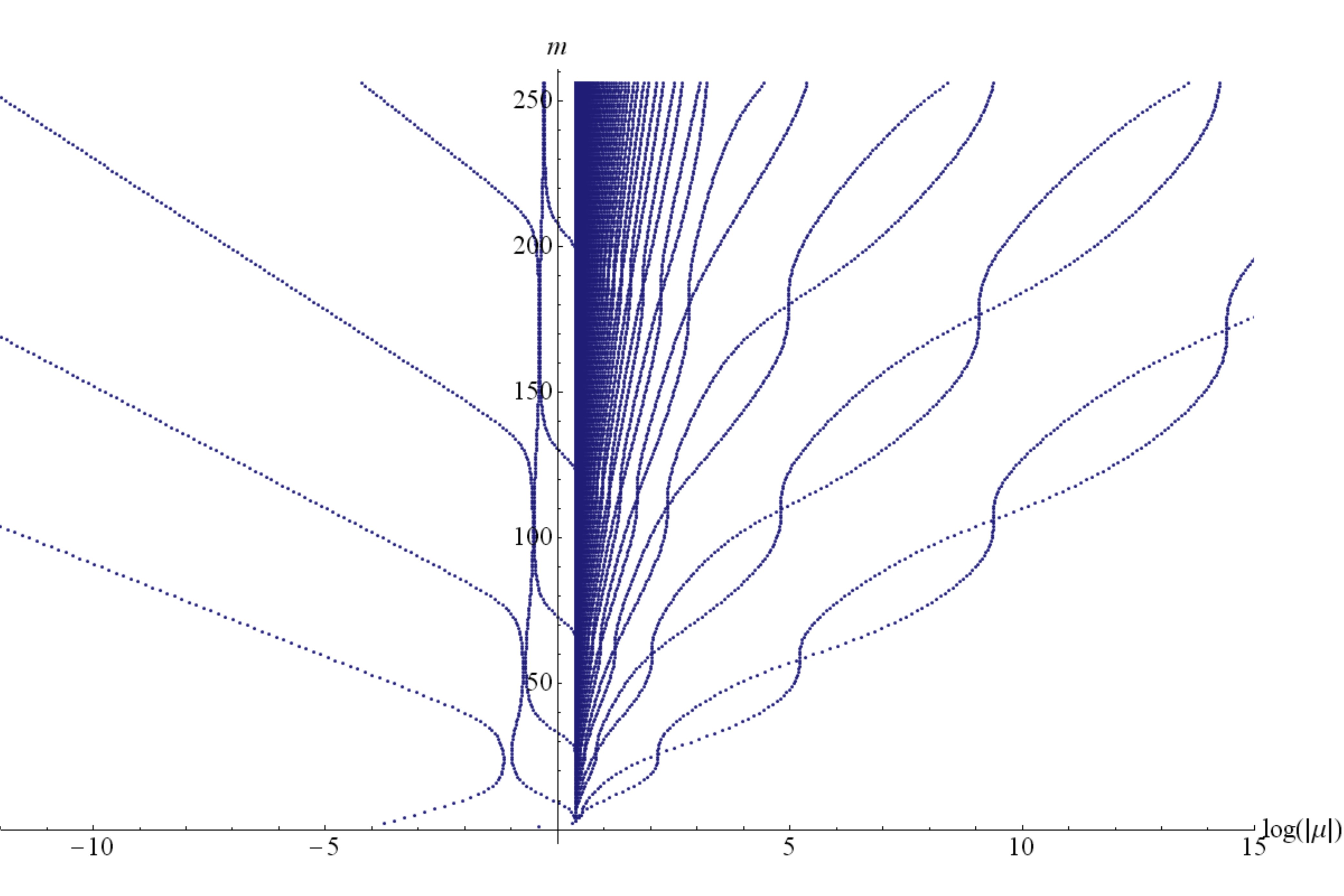}\\
\parbox{.94\textwidth}{\caption{$\lamS_{2,m}(\zetastar),\ m=1,\dots,256.$}
\label{picm2}}
\end{figure}
\begin{figure}[phtb]
\includegraphics[width=1\textwidth]{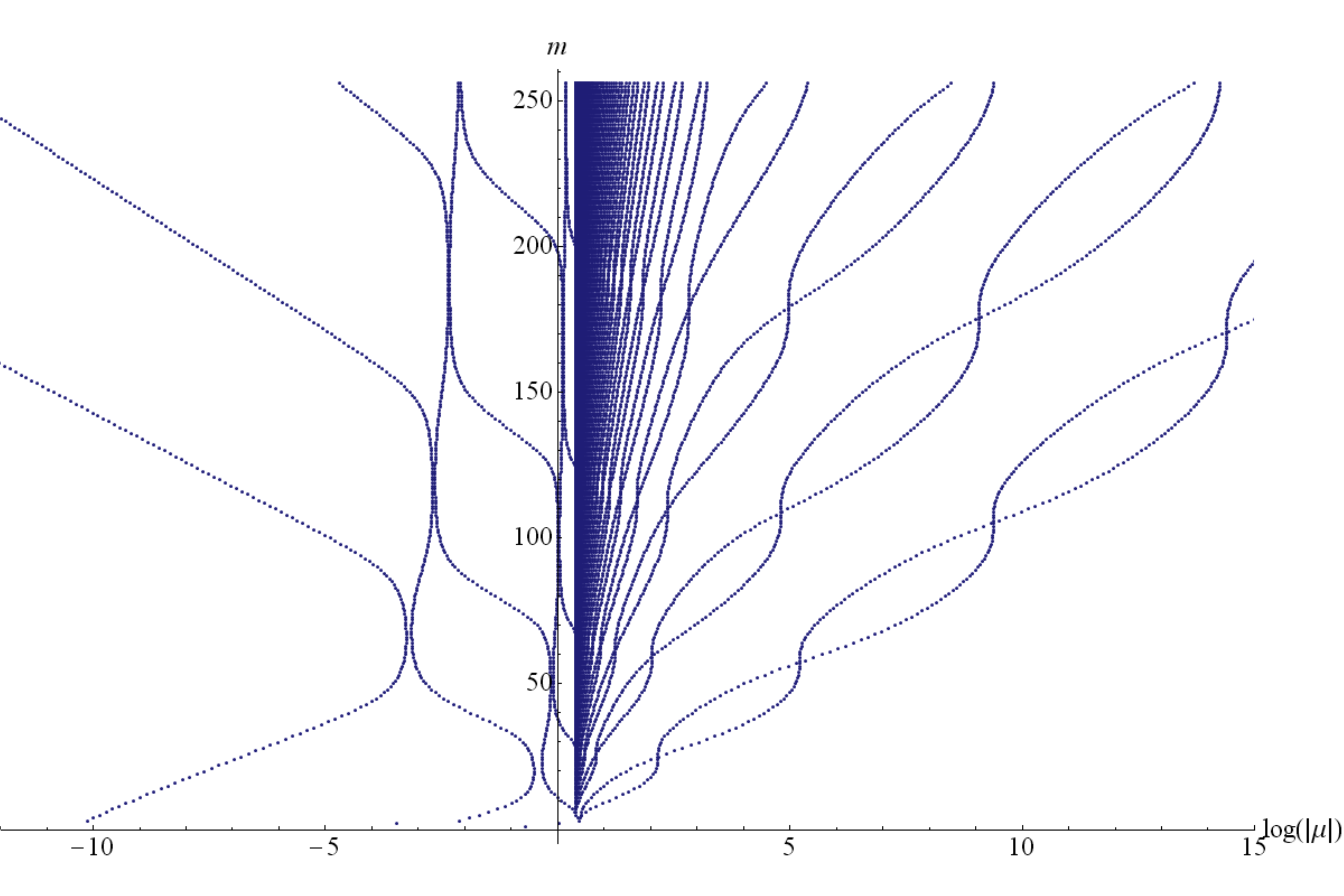}\\
\parbox{.94\textwidth}{\caption{$\lamS_{3,m}(\zetastar),\ m=1,\dots,256.$}
\label{picm3}}
\end{figure}
\begin{figure}[phtb]
\includegraphics[width=1\textwidth]{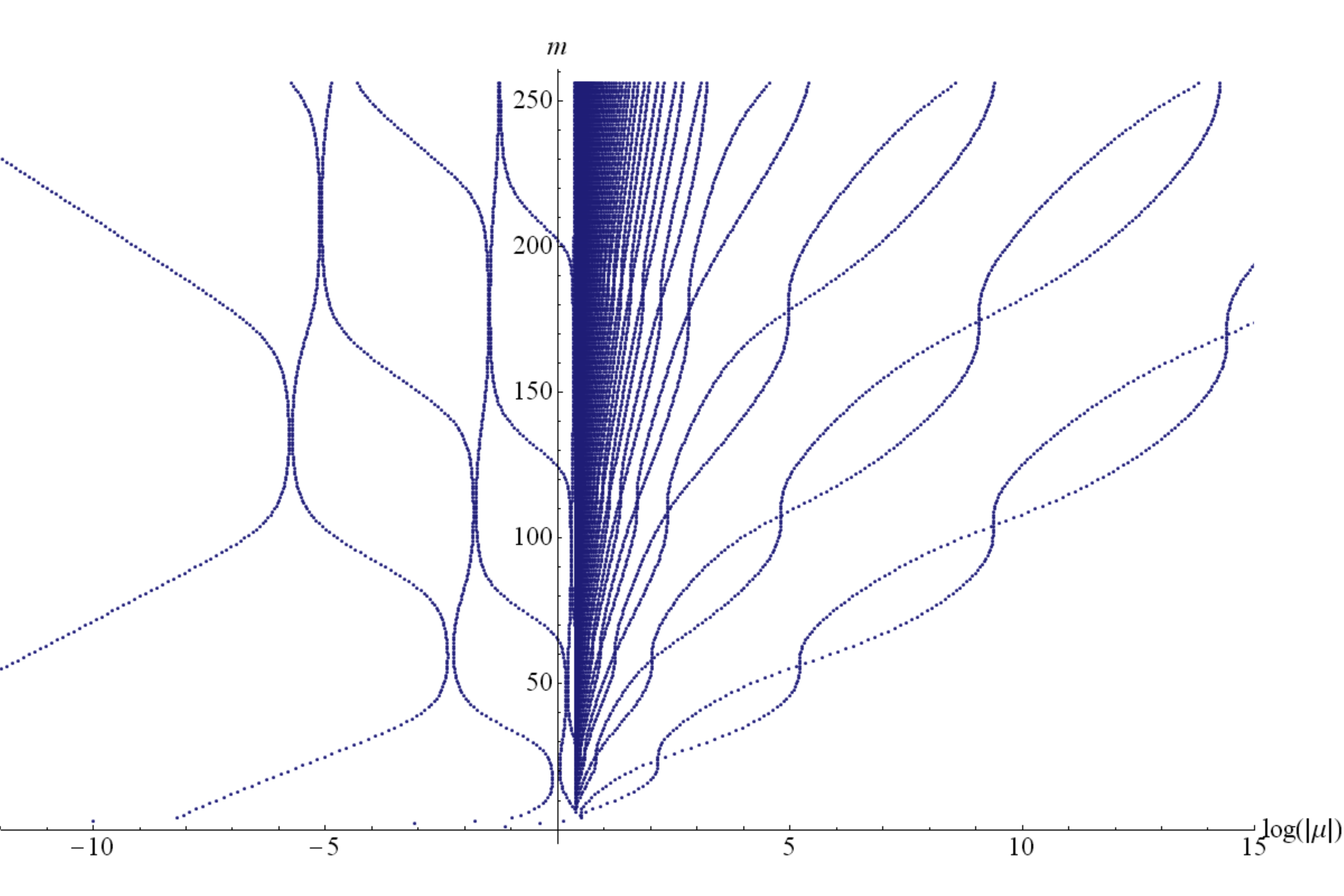}\\
\parbox{.94\textwidth}{\caption{$\lamS_{4,m}(\zetastar),\ m=1,\dots,256.$}
\label{picm4}}
\end{figure}
\begin{figure}[phtb]
\includegraphics[width=1\textwidth]{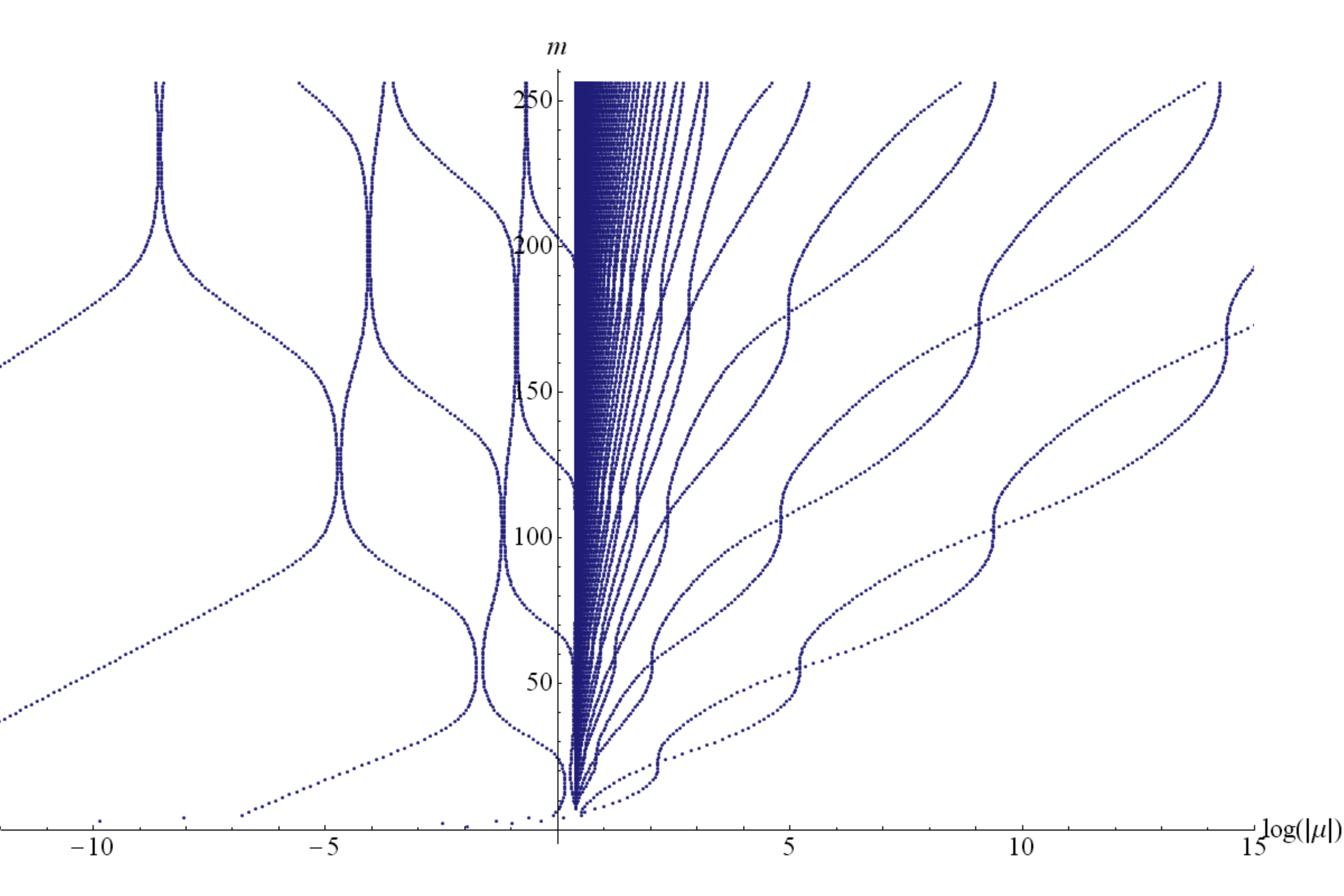}\\
\parbox{.94\textwidth}{\caption{$\lamS_{5,m}(\zetastar),\ m=1,\dots,256.$}
\label{picm5}}
\end{figure}
\begin{figure}[phtb]
\includegraphics[width=1\textwidth]{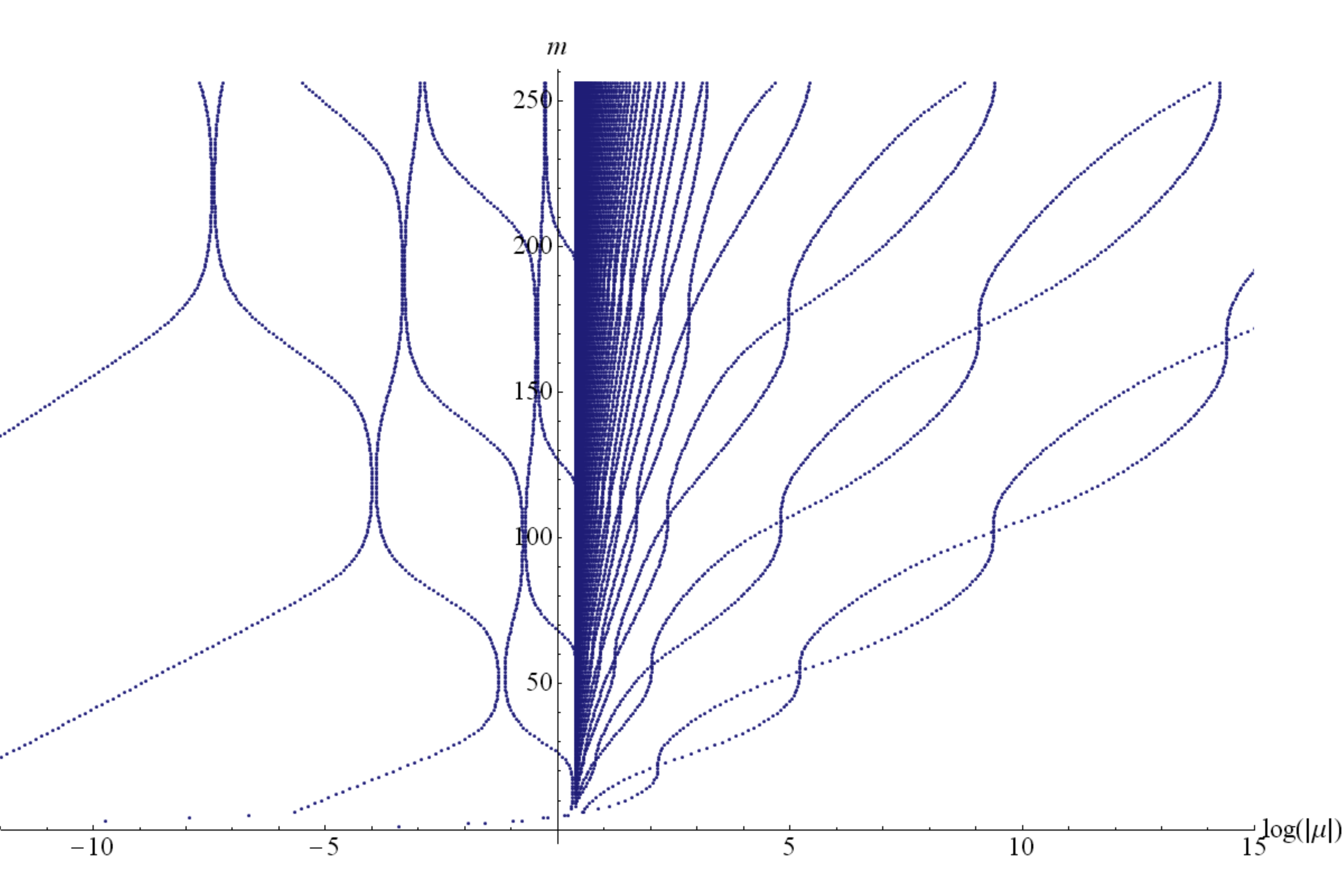}\\
\parbox{.94\textwidth}{\caption{$\lamS_{6,m}(\zetastar),\ m=1,\dots,256.$}
\label{picm6}}
\end{figure}

Figures \ref{picm1}--\ref{picm6} show
spectra $\lamS_{1,m}(\zetastar)$,
\dots,
$\lamS_{2,m}(\zetastar)$ for $m=1,\dots,256$
(higher resolution version of these pictures
can be downloaded from \cite{hiddensite} as
well as some animations showing these
spectra
 and thus
revealing another kind of  ``hidden life of
Riemann's zeta function'').

The pictures and the animations show that
with the growth of $m$ some elements of
$\lamS_{l,m}(\zetastar)$ go to $-\infty$
while  others go to $+\infty$; the former
will be called \emph{electrons} and the
latter will be called \emph{trains} (we
postpone formal definition of splitting
$\lamS_{l,m}(\zetastar)$ into \emph{lower
part} $\lamSlow_{l,m}(\zetastar)$,
consisting of the electrons,
 and \emph{upper part} $\lamSup_{l,m}(\zetastar)$ consisting of
the trains).

The names ``electrons'' and ``trains'' were
suggested by the following visual patterns.
The electrons behave like charged particles,
namely, they bounce. The trains all go in
pairs (a surprising feature!) and every now
and then they overtake one another.

\subsection{New Conjectures}

The above pictures suggest the following
conjectures.

\

{\bf Conjecture 2A.} \emph{For all $l$}
\begin{equation}
\lim_{m\rightarrow\infty} \max
(\lamS_{l,m}(\zetastar))=+\infty.
\end{equation}

\

{\bf Conjecture 2B.} \emph{For all $l$}
\begin{equation}
\lim_{m\rightarrow\infty} \min
(\lamS_{l,m}(\zetastar))=-\infty.
\end{equation}

\

It is impossible to see from the above
pictures whether for the $\mu$-spectra there
is a counterpart of Conjecture 1F about the
$\lambda$-spectra. To make this clearer, in
analogy with this conjecture, let us assign
to each point of $\lamS_{l,m}(\zetastar)$
the weight $\frac{1}{m}$, and denote by
$\mu^{\zetastar}_{l,m}(x)$ the corresponding
discrete measure on real numbers. Further,
let $F^{\zetastar}_{l,m}(x)$ denote the
corresponding distribution function. In
terms of these functions we have

\

{\bf  $\mathbf{RH_{\emph{l}}^w}$ (version
6)}.
\begin{equation}
\lim_{m\rightarrow \infty} \left(
\int_{-\infty}^{+\infty}x\, \myd
F^{\zetastar}_{l,m}(x)\right)=\log(\W_l).
\label{vers6}\end{equation}

\begin{figure}[phtb]
\includegraphics[width=.47\textwidth]{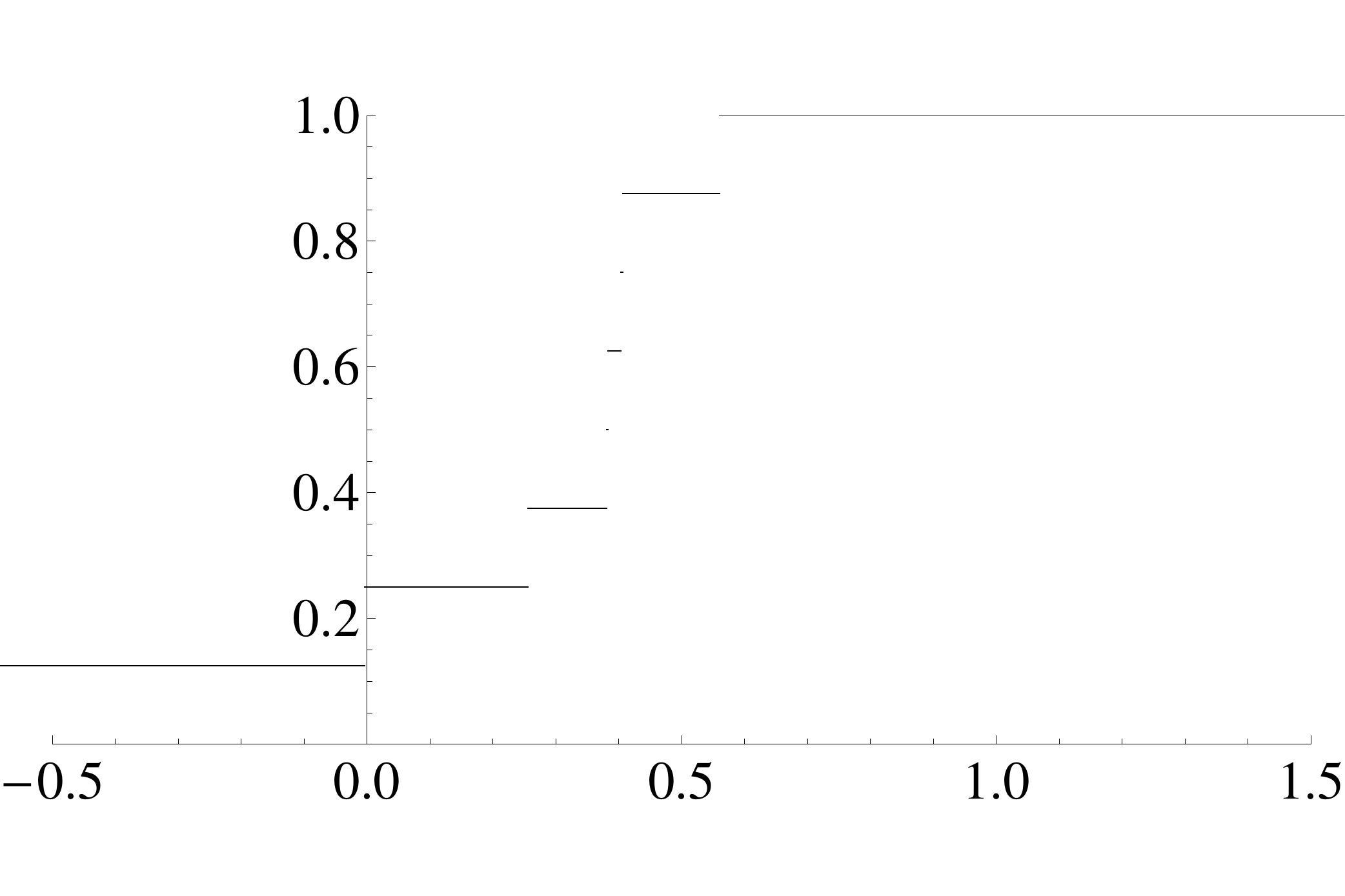}\hfill
\includegraphics[width=.47\textwidth]{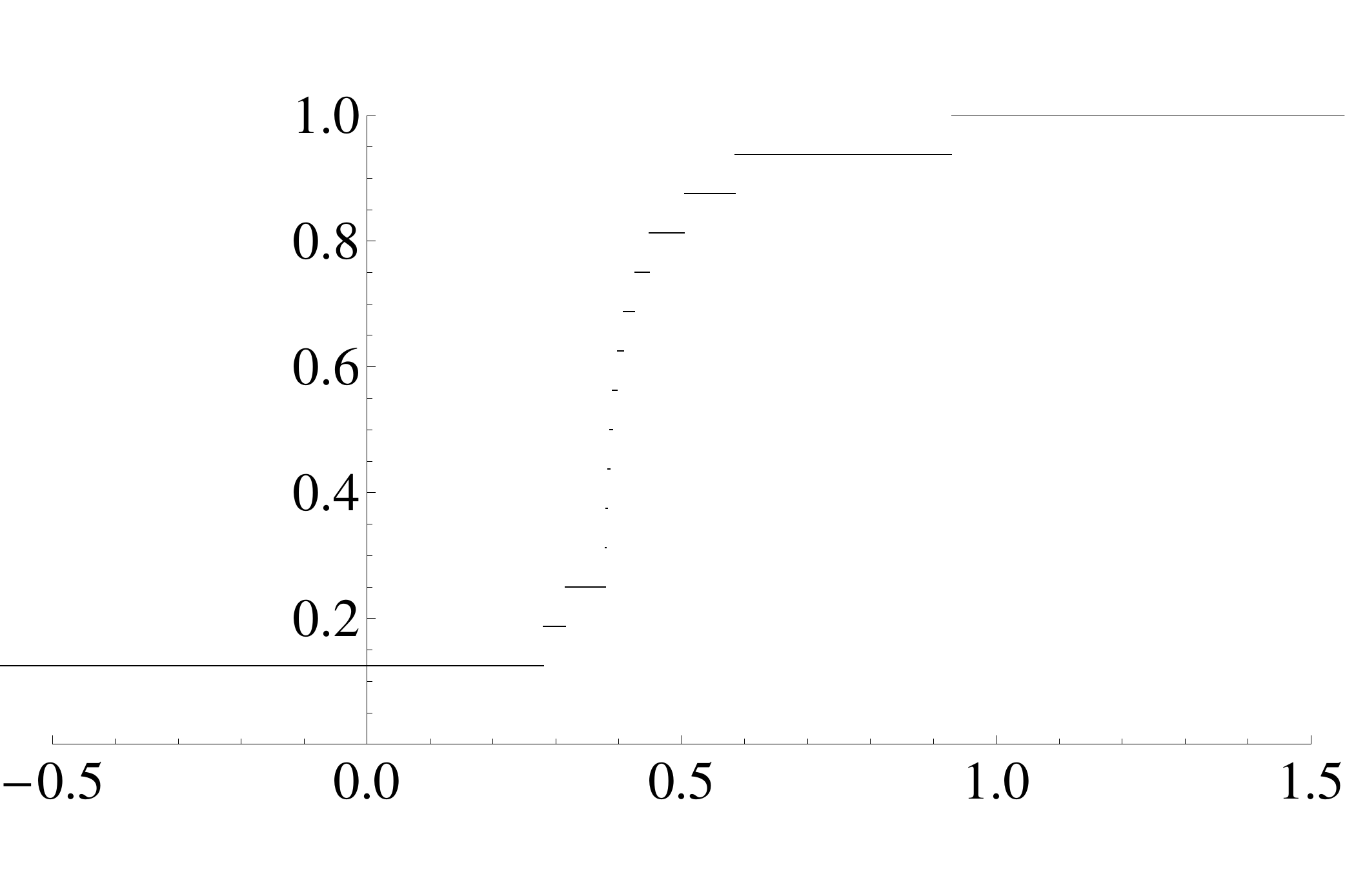}\\
\parbox{.47\textwidth}{\caption{$F^{\zetastar}_{1,8}(x).$}\label{picdist18}}\hfill
\parbox{.47\textwidth}{\caption{$F^{\zetastar}_{1,16}(x).$}\label{picdist116}}
\end{figure}
\begin{figure}[phtb]
\includegraphics[width=.47\textwidth]{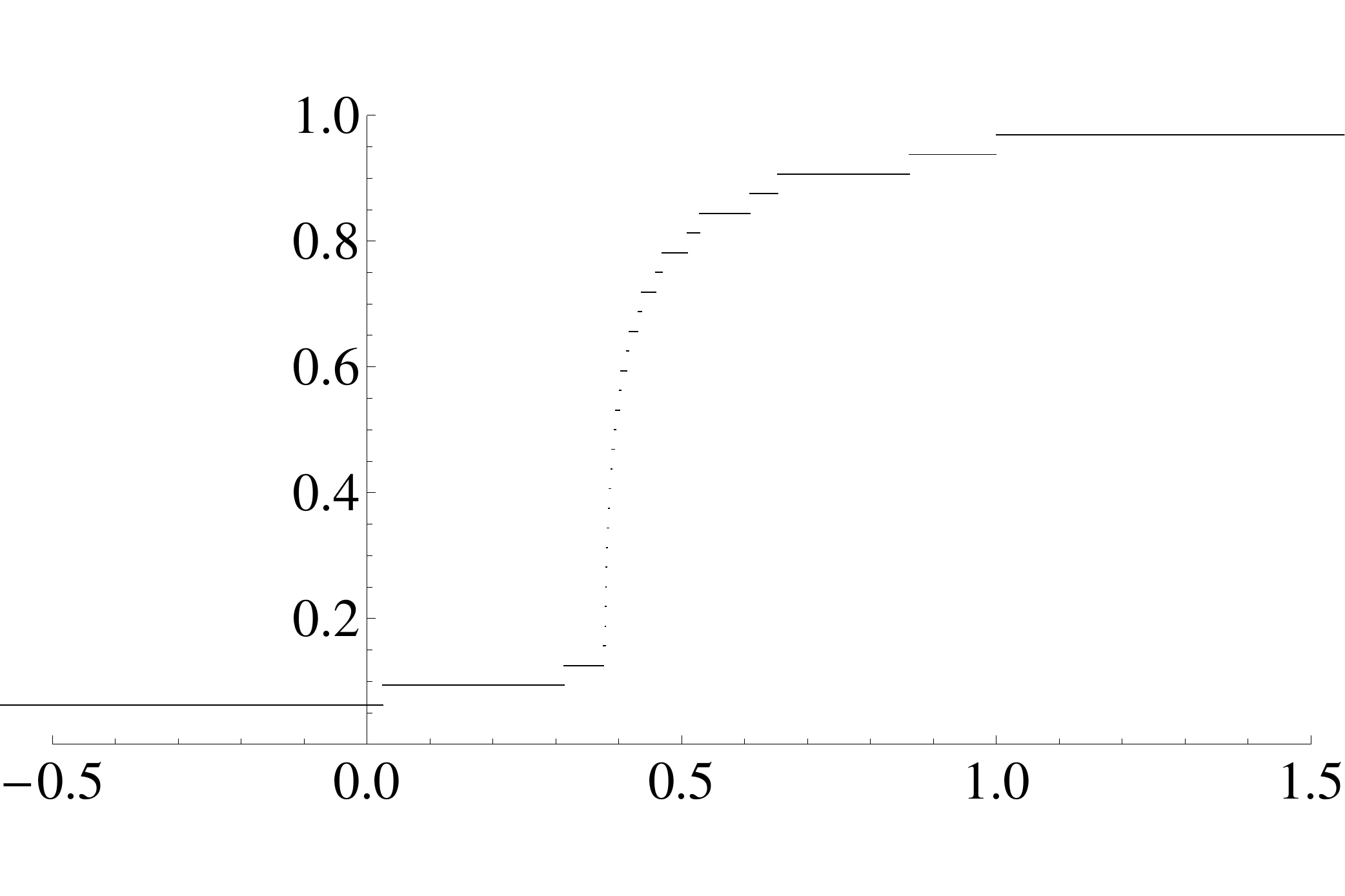}\hfill
\includegraphics[width=.47\textwidth]{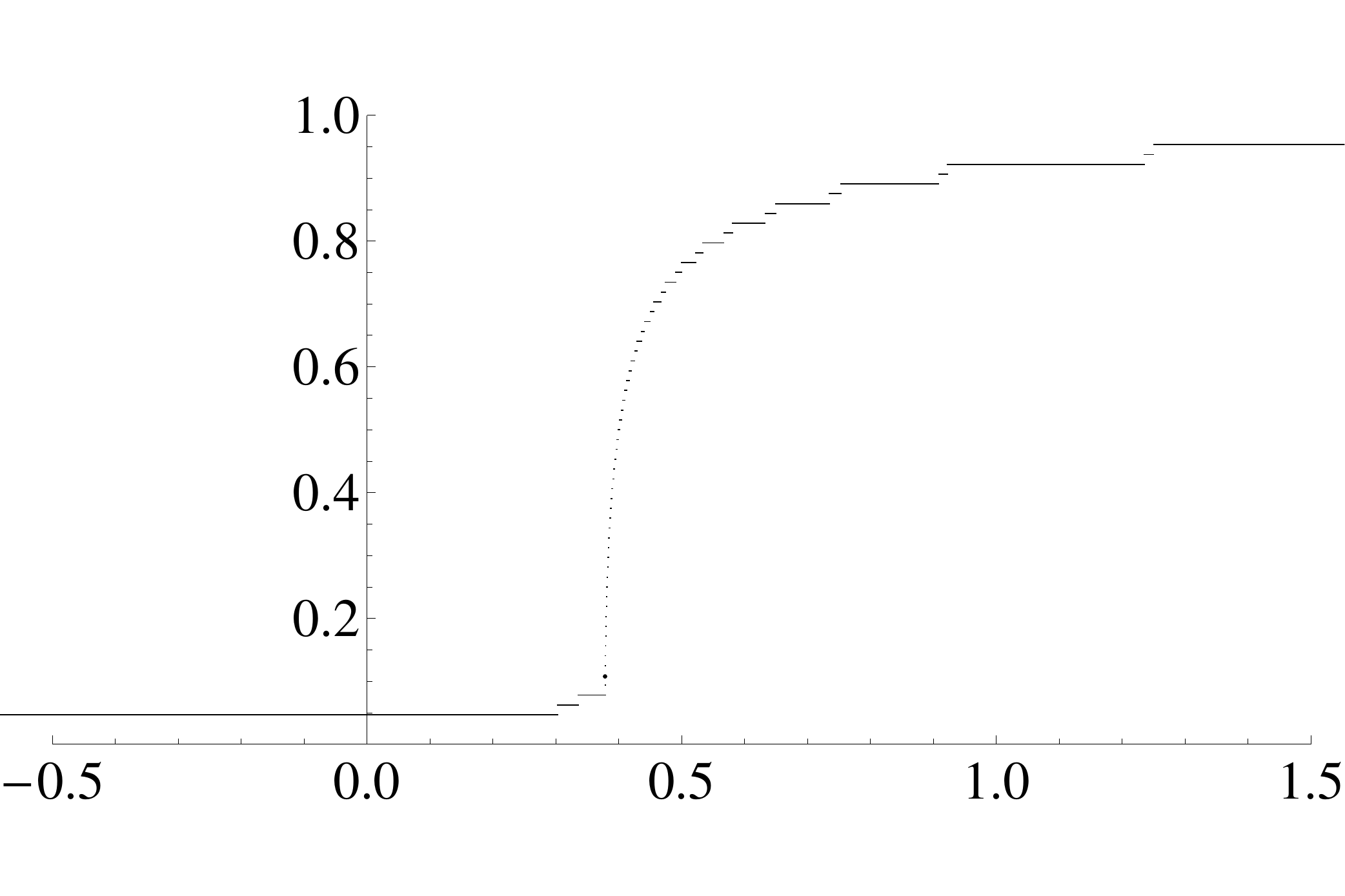}\\
\parbox{.47\textwidth}{\caption{$F^{\zetastar}_{1,32}(x).$}\label{picdist132}}\hfill
\parbox{.47\textwidth}{\caption{$F^{\zetastar}_{1,64}(x).$}\label{picdist164}}
\end{figure}
\begin{figure}[phtb]
\includegraphics[width=.47\textwidth]{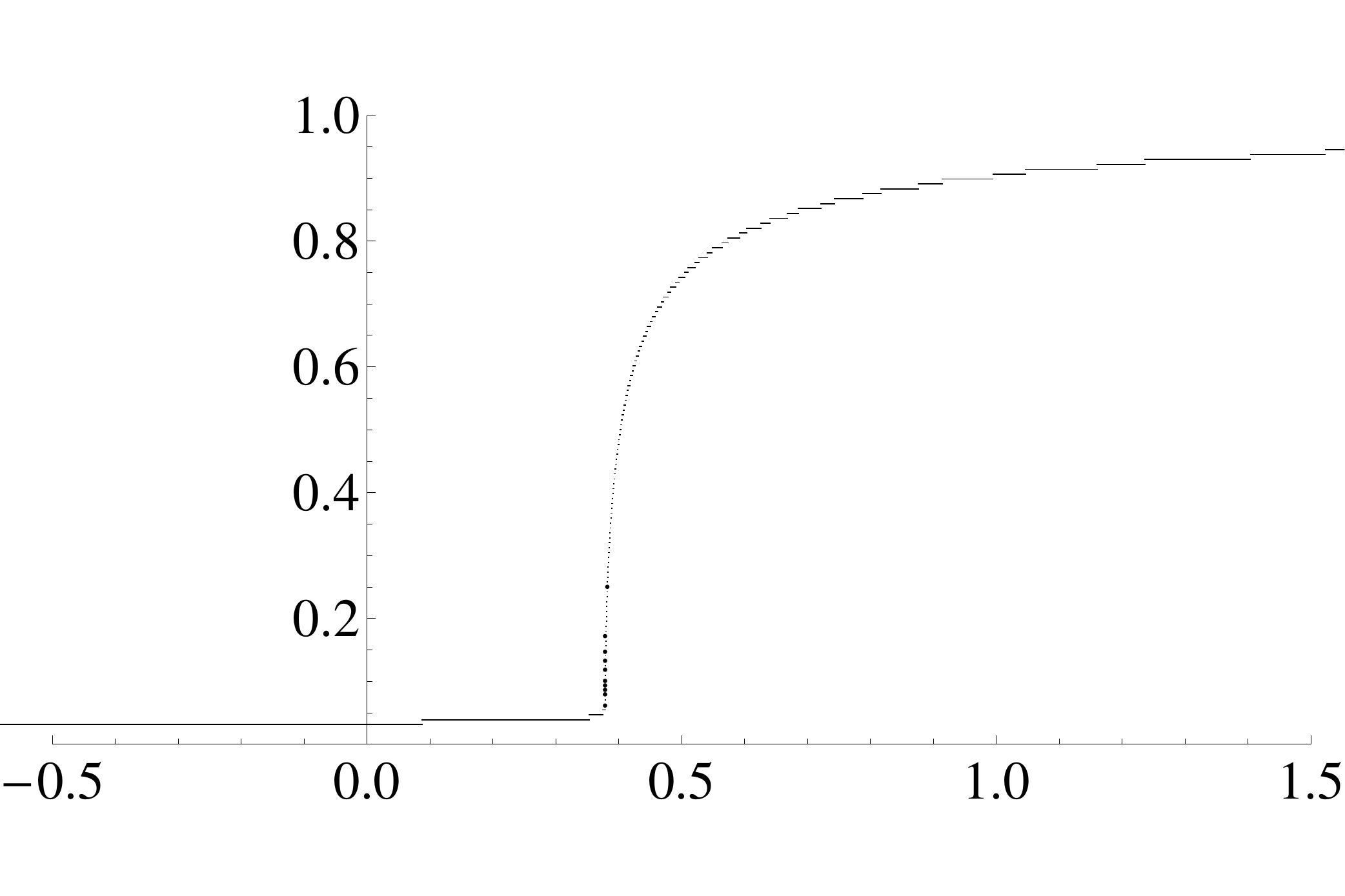}\hfill
\includegraphics[width=.47\textwidth]{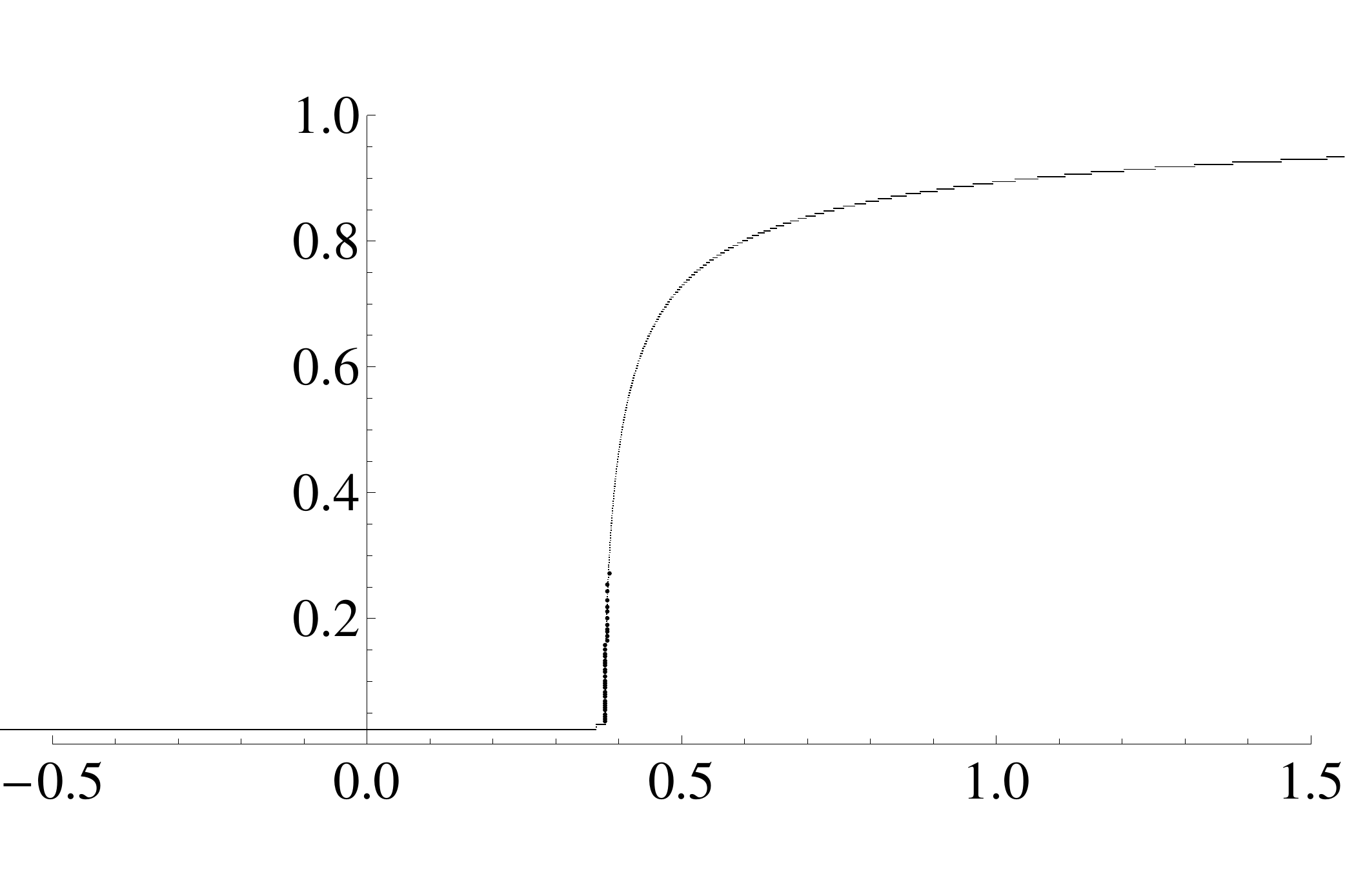}\\
\parbox{.47\textwidth}{\caption{$F^{\zetastar}_{1,128}(x).$}\label{picdist1128}}\hfill
\parbox{.47\textwidth}{\caption{$F^{\zetastar}_{1,256}(x).$}\label{picdist1256}}
\end{figure}
\begin{figure}[phtb]
\includegraphics[width=.47\textwidth]{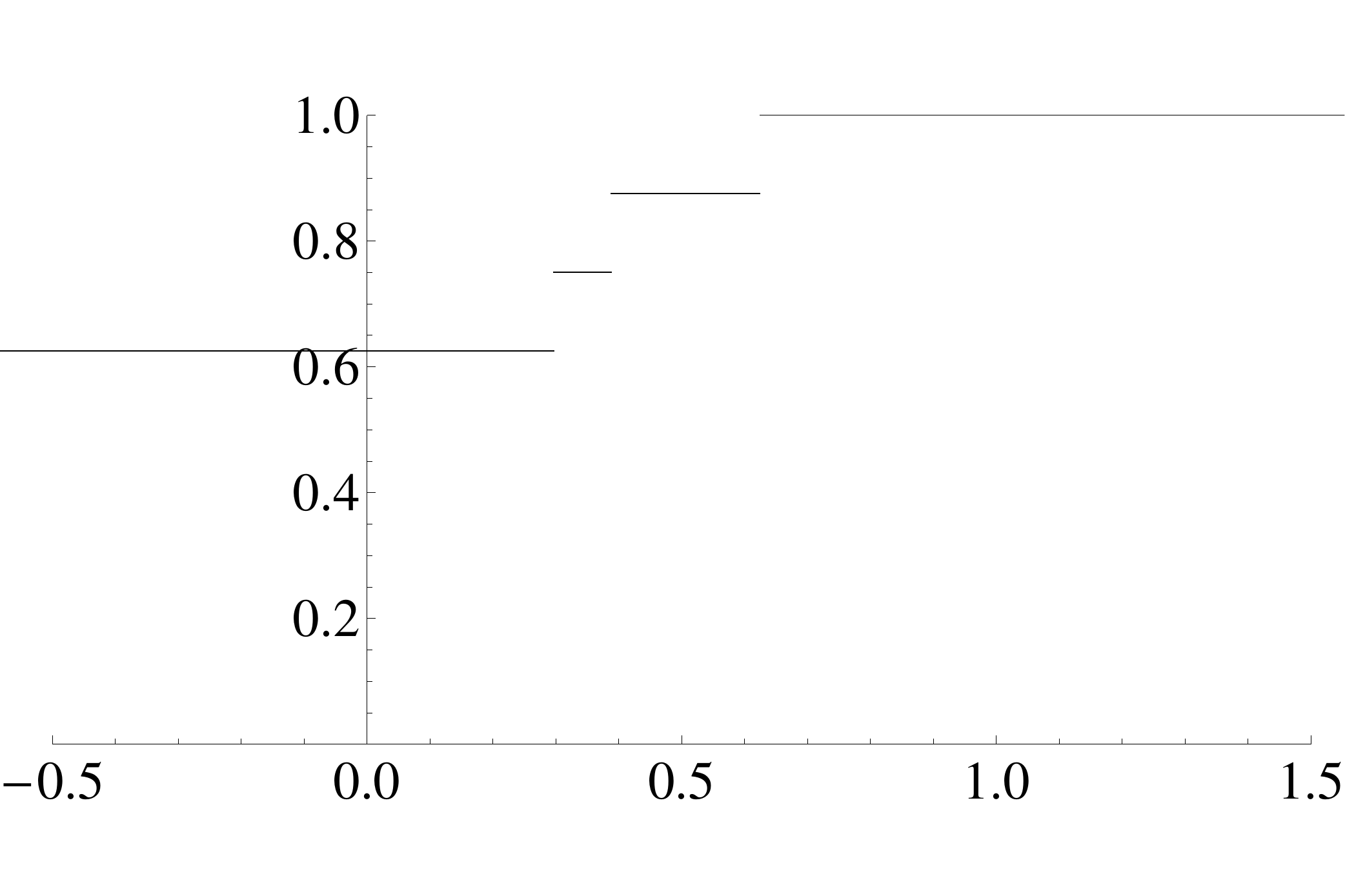}\hfill
\includegraphics[width=.47\textwidth]{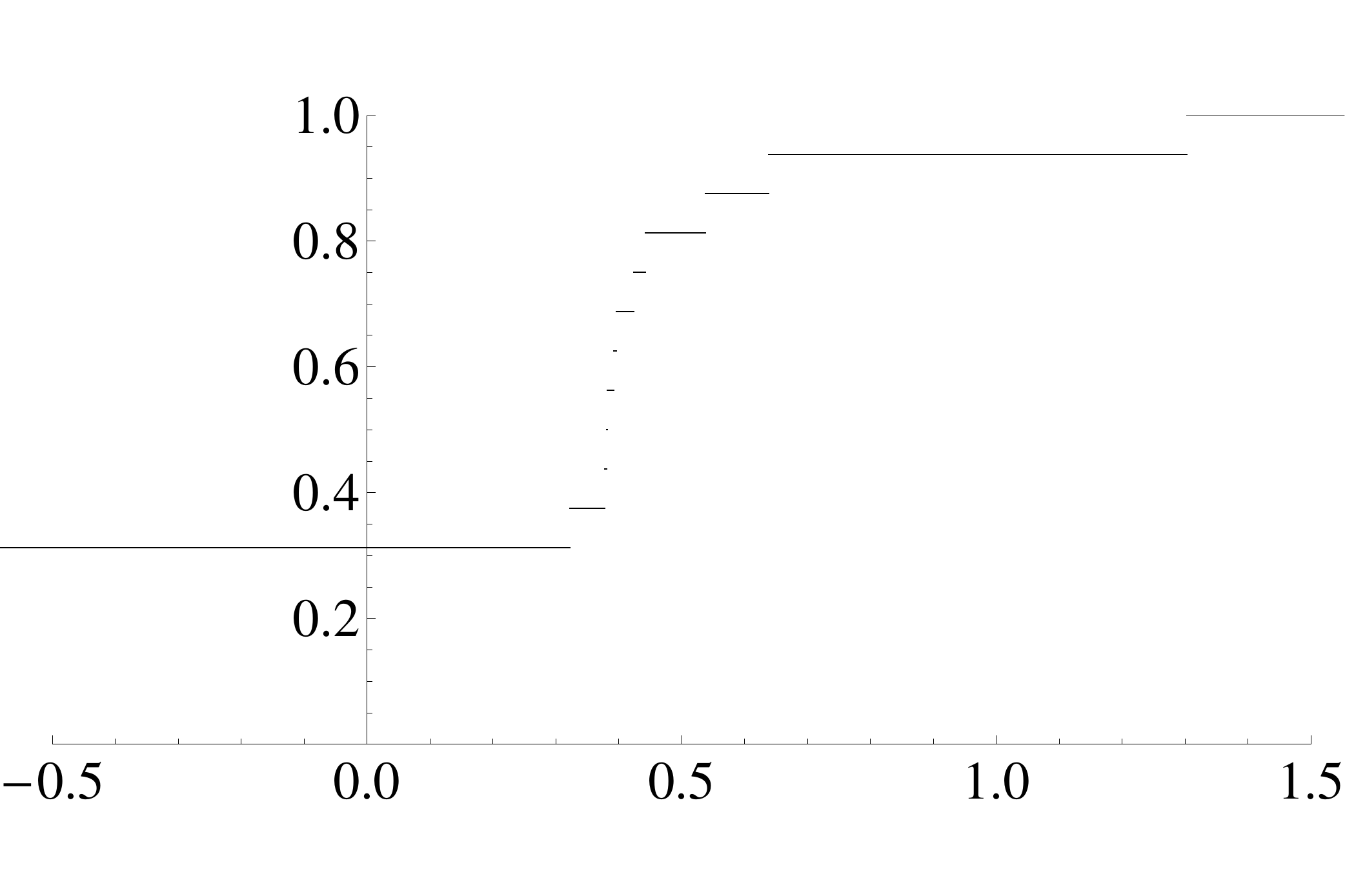}\\
\parbox{.47\textwidth}{\caption{$F^{\zetastar}_{6,8}(x).$}\label{picdist68}}\hfill
\parbox{.47\textwidth}{\caption{$F^{\zetastar}_{6,16}(x).$}\label{picdist616}}
\end{figure}
\begin{figure}[phtb]
\includegraphics[width=.47\textwidth]{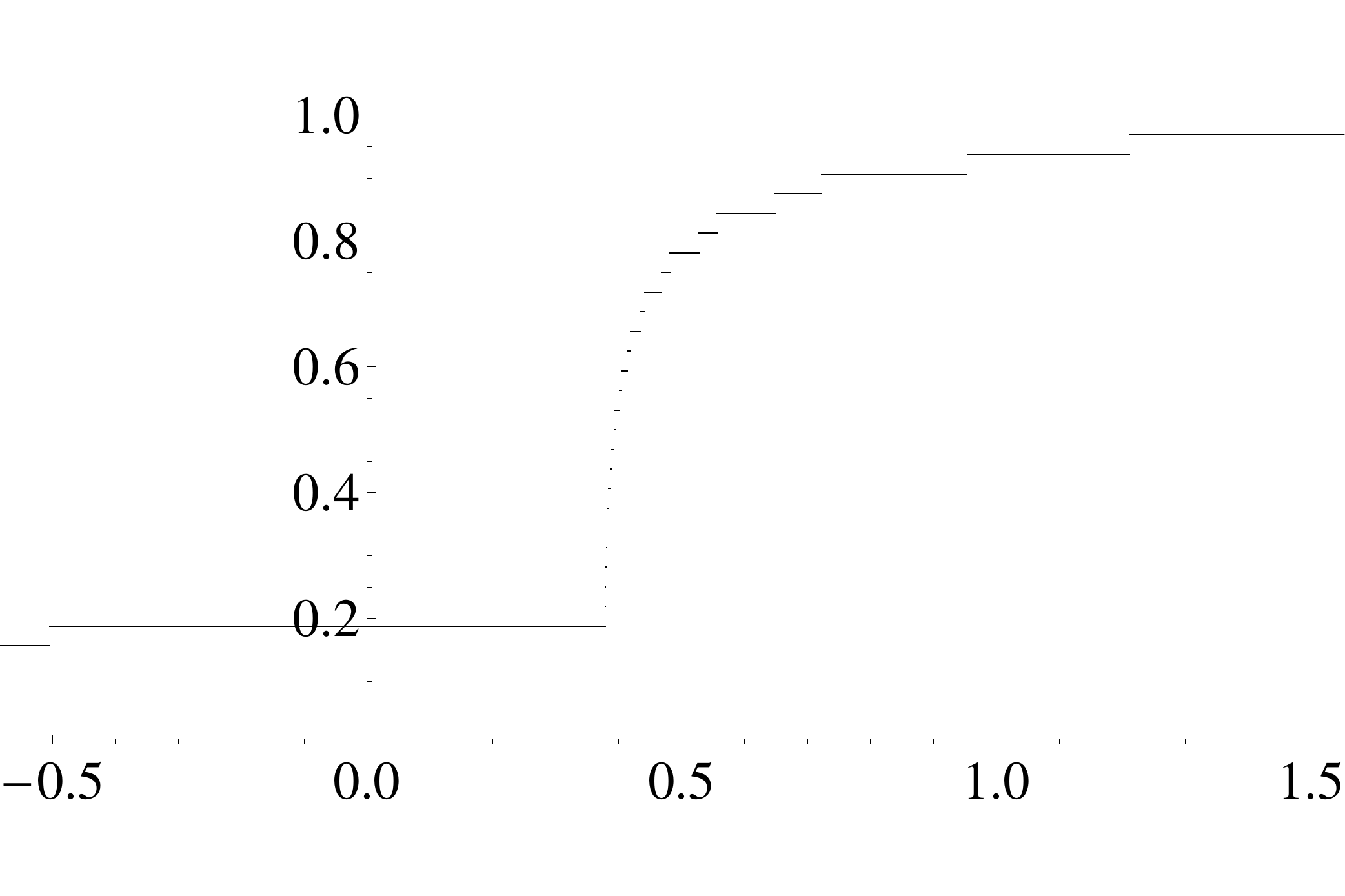}\hfill
\includegraphics[width=.47\textwidth]{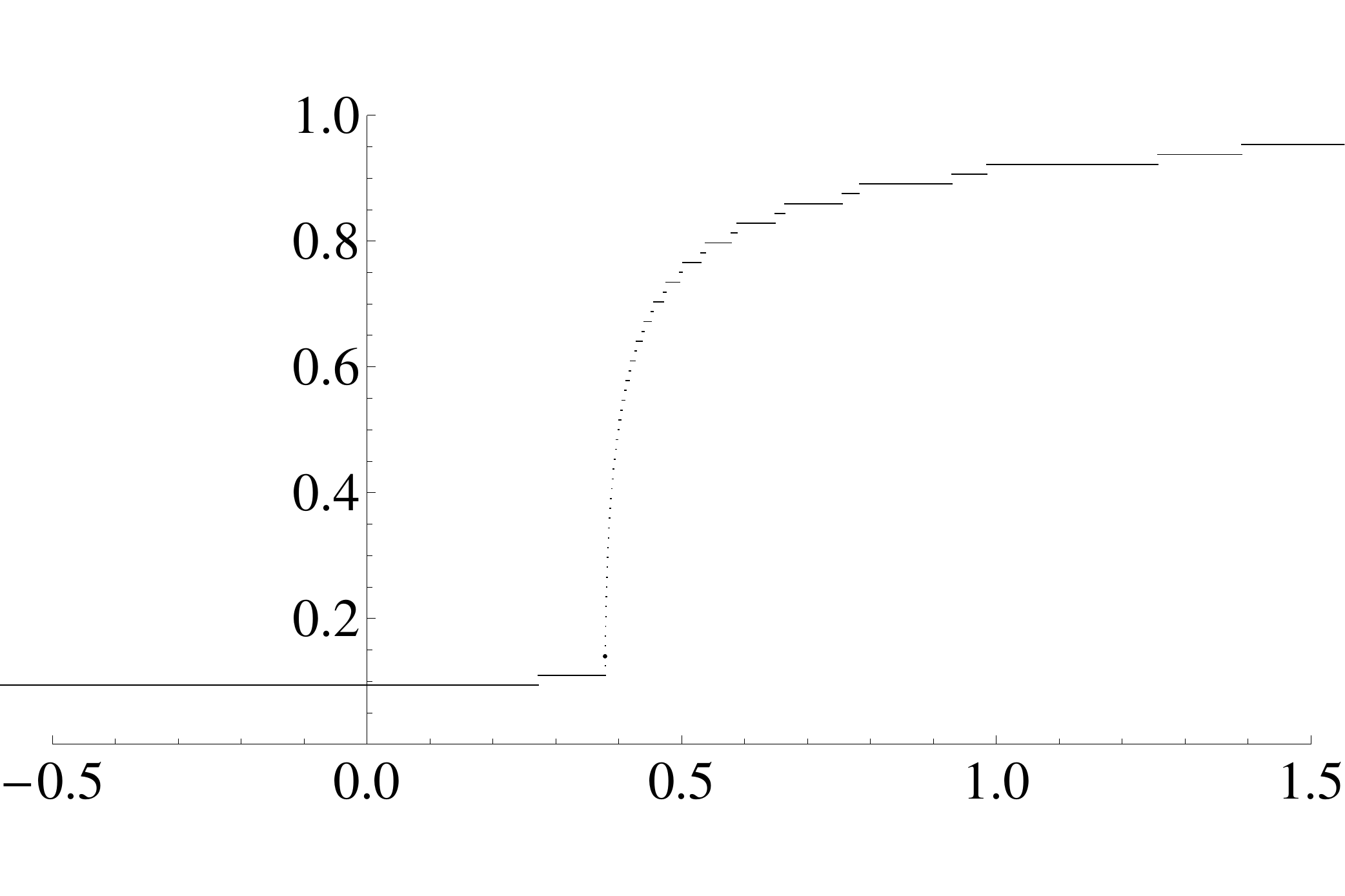}\\
\parbox{.47\textwidth}{\caption{$F^{\zetastar}_{6,32}(x).$}\label{picdist632}}\hfill
\parbox{.47\textwidth}{\caption{$F^{\zetastar}_{6,64}(x).$}\label{picdist664}}
\end{figure}
\begin{figure}[phtb]
\includegraphics[width=.47\textwidth]{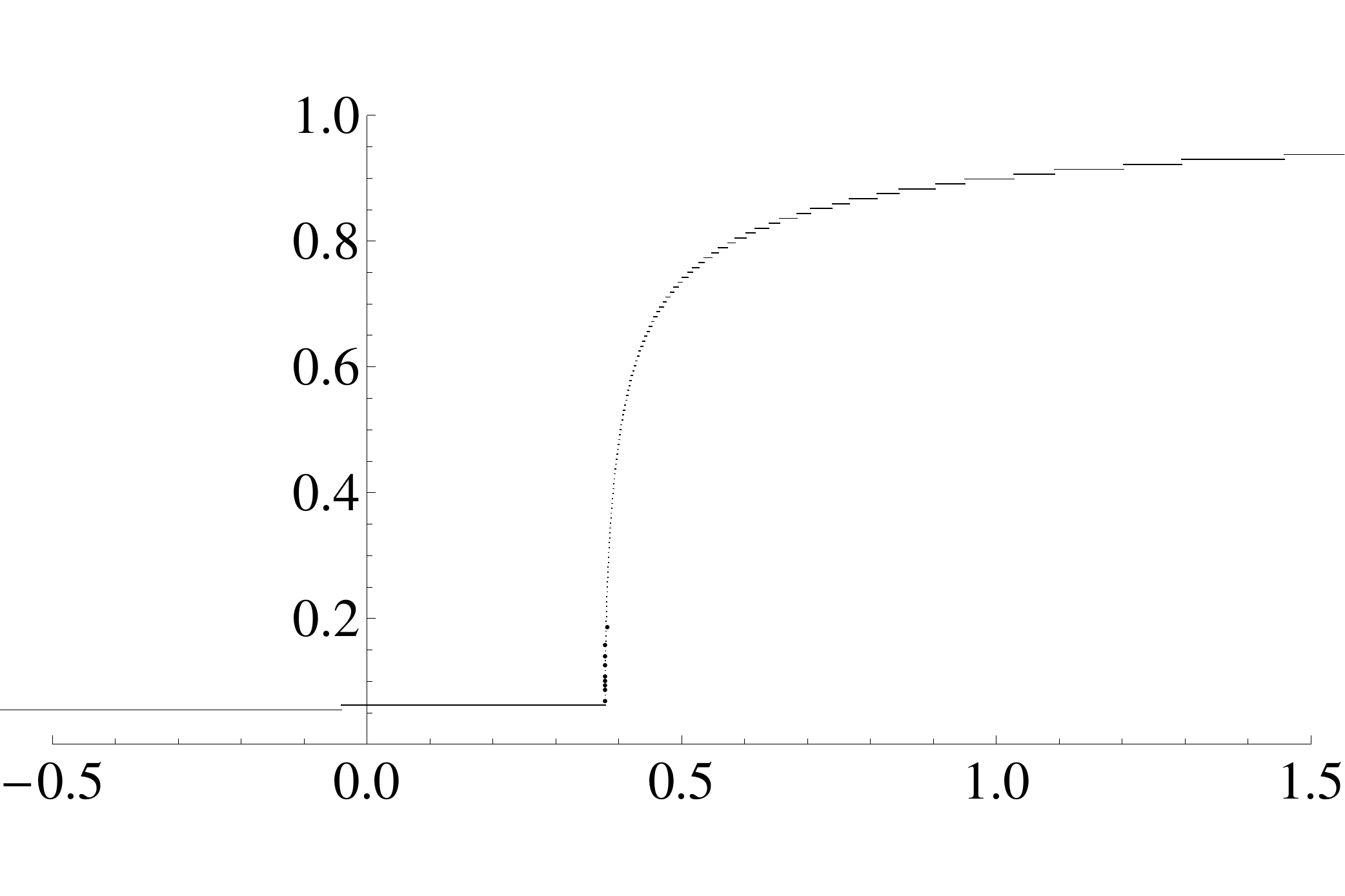}\hfill
\includegraphics[width=.47\textwidth]{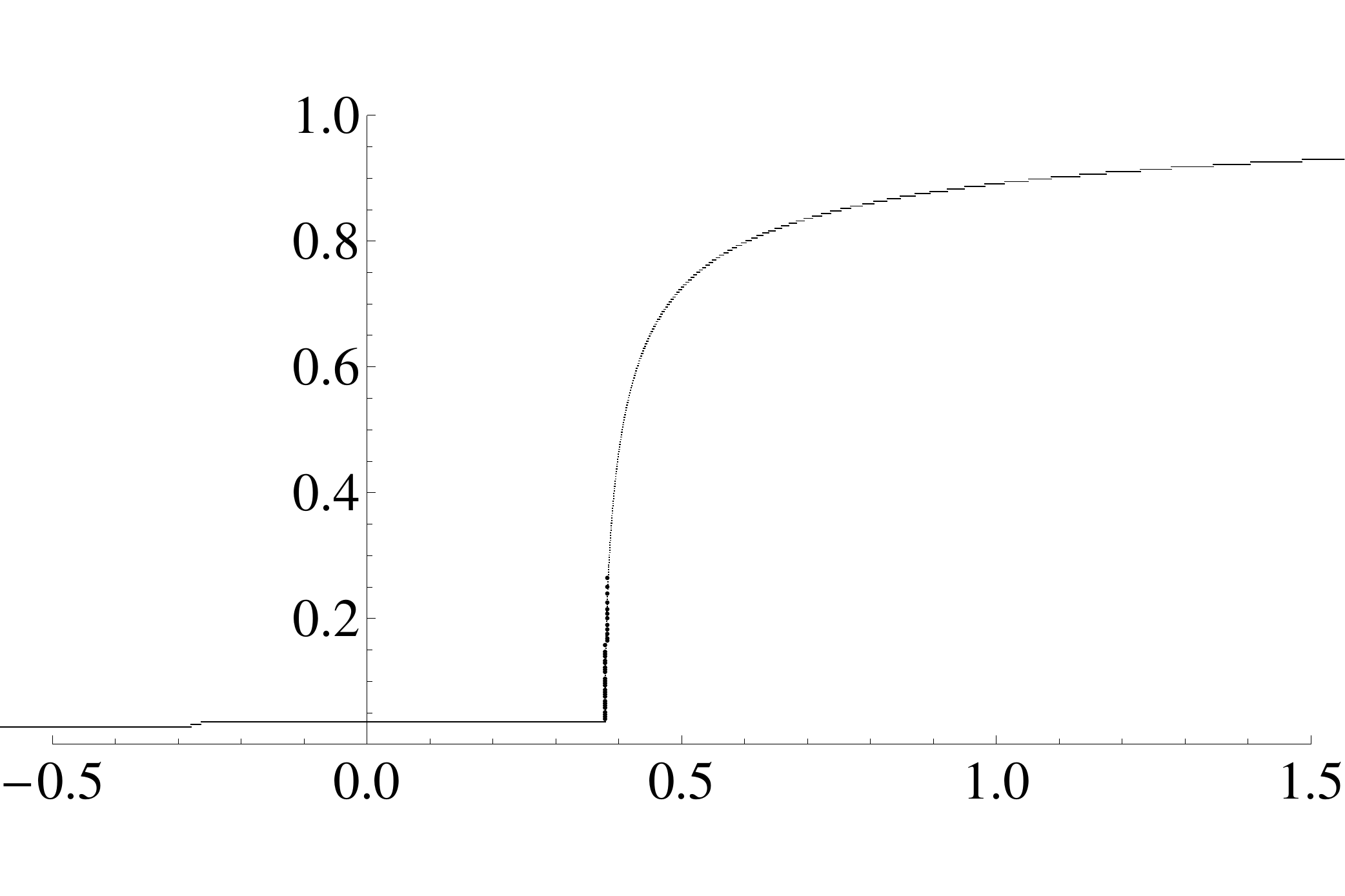}\\
\parbox{.47\textwidth}{\caption{$F^{\zetastar}_{6,128}(x).$}\label{picdist6128}}\hfill
\parbox{.47\textwidth}{\caption{$F^{\zetastar}_{6,256}(x).$}\label{picdist6256}}
\end{figure}

Figures  \ref{picdist18}--\ref{picdist6256} show
these functions for $l=1,6$ and  $m=8,16,32,64,256$.
These pictures suggest

\

{\bf Conjecture 2C.}
 \emph{For every $l$ functions
$F^{\zetastar}_{l,m}(x)$  have, as
$m\rightarrow \infty$, the pointwise
limiting continuous distribution function
$F_{l}^{\zetastar}\!(x)$.}

\

This  is an analog of part $\mathrm{1F'}$ of
Conjecture 1F for $\lambda$-spectra.
However, it seems that part  $\mathrm{1F''}$
of this conjecture has no analog for
$\mu$-spectra. According to \eqref{vers6},
such an analog would say that for all $l$
\begin{equation}
\int_{-\infty}^{+\infty}x\, \myd
F^{\zetastar}_{l}(x)=\log(\W_l).
\label{wrongconj}\end{equation} But theses
integrals do not seem  to exist:

\

{\bf Conjecture 2D.} \emph{For every $l$}
\begin{eqnarray}
\int_{-\infty}^{0}x\, \myd
F^{\zetastar}_{l}(x)=-\infty ,&&
\int_{0}^{+\infty}x\, \myd
F^{\zetastar}_{l}(x)=+\infty.
\label{intinfty}
\end{eqnarray}

\

This implies that the validity of
\eqref{vers6} should be due to some fine
correlation between eigenvalues from
$\lamSlow_{l,m}(\zetastar)$ and
$\lamSup_{l,m}(\zetastar)$. The subtleness
of this correlation follows from (very surprising)

\

{\bf Conjecture 2E.} \emph{All distribution
functions  $F^{\zetastar}_{l}(x)$ coincide;
that is, for all~$l$ and $x$
\begin{equation}
F^{\zetastar}_{l}(x)=F^{\zetastar}(x)
\end{equation}
for some continuous  distribution function
$F^{\zetastar}(x)$.}

\section*{Acknowledgement}

\

The author is very grateful to Martin Davis
for some help with the English.

\end{document}